\documentclass[11pt,dvips,epsfig,reqno,graphics]{amsart}
\usepackage{a4wide,latexsym,amssymb,amsthm,amsmath,amsfonts,url}

\input{tcilatex}

\begin{document}
\title[Semi-invariant submanifolds in metric geometry of affinors]{%
Semi-invariant submanifolds in metric geometry of affinors}
\author[N.-C. Chiriac]{Novac-Claudiu Chiriac $^\dag$}
\thanks{$^\dag$ Corresponding Author}
\address[N.-C. Chiriac]{Department of Mathematics \\
University Constantin Br\^ancusi \\
Bd. Republicii, Nr. 1 \\
T\^argu-Jiu, 210152 \\
Romania}
\email{novac (at) utgjiu.ro}
\author[M.~Crasmareanu]{Mircea Crasmareanu}
\address[M.~Crasmareanu]{Faculty of Mathematics \\
'Al. I. Cuza' University of Iasi\\
Bd. Carol I, nr. 11\\
Iasi, 700506 \\
Romania, \quad \url{http://www.math.uaic.ro/~mcrasm}}
\email{mcrasm (at) uaic.ro}
\date{August 19, 2011}

\begin{abstract}
We introduce a generalization of structured manifolds as the most general
Riemannian metric $g$ associated to an affinor (tensor field of $(1, 1)$%
-type) $F$ and initiate a study of their semi-invariant submanifolds. These
submanifolds are generalizations of CR-submanifolds of almost complex
geometry and semi-invariant submanifolds of several in\-te\-res\-ting
geometries (almost product, almost contact and others). We characterize the
integrability of both invariant and anti-invariant distribution; the special
case when $F$ is covariant constant with respect to $g$ gives major
simplifications in computations.
\end{abstract}

\maketitle



\section*{Introduction}

The geometry of manifolds endowed with geometrical structures has been
intensively stu\-died and several important results have been published, see
Yano-Kon \cite{14}. The more important classes of such manifolds are formed
by almost complex, almost product, almost contact, almost paracontact
manifolds for which the cited book offers a good introduction. The geometry
of submanifolds in these manifolds is very rich and interesting, as well,
see for example the classical \cite{7} or the more recent survey \cite{8}.
CR-submanifolds introduced by Bejancu in \cite{2} (for almost complex
geometry) respectively \cite{5} (for almost contact geometry) had a great
impact on the developing of the theory of submanifolds in these ambient
manifolds; a proof of this fact is given by the books \cite{4} and \cite{13}.

\medskip

In the present paper we first introduce the concept of $(g,F,\mu )$-manifold
which contains as particular cases all the above types of structures. Then
we study semi-invariant sub\-ma\-ni\-folds of a $(g,F,\mu )$-manifold, which
are extensions of CR-submanifolds to this general class of manifolds. We
find necessary and sufficient conditions for the integrability of both
distributions on a semi-invariant submanifold, see Theorems 3.1 and 3.3.
Particularly, we prove that some semi-invariant submanifolds carry a natural
foliation, in Theorem 4.4 and we obtain characterizations of totally
geodesic foliations on semi-invariant submanifolds in Theorems 4.8 and 4.10.
For a particular value of the real parameter $\mu $ we can connect our study
with the almost symplectic geometry and this fact opens some possible
further applications in physical sciences having as an example the
relationship between CR-structures and Relativity pointed out in the last
Chapter of \cite{4}.

\medskip

This work is dedicated to Professor Aurel Bejancu on the occasion of his
65th birthday. His ideas represent a starting point for several important
studies as the present Bibliography partially shows.

\section{Metric geometry of affinors and submanifolds}

Let $M$ be an $m$-dimensional manifold for which we denote by $C^{\infty
}(M) $ the algebra of smooth functions on $M$ and by $\Gamma (TM)$ the $%
C^{\infty }(M)$-module of smooth sections of the tangent bundle $TM$ of $M$;
let $X,Y,Z,...$ denote such vector fields. We use the same notation $\Gamma
(V)$ for any other vector bundle $V$ over $M$. Let also ${\mathcal{T}}%
_{1}^{1}M$ be the $C^{\infty }(M)$-module of $\Gamma (TM\otimes T^{\ast }M)$
i.e. the real space of tensor fields of $(1,1)$-type on $M$. Let consider a
fixed $F\in {\mathcal{T}}_{1}^{1}M$ usually called \textit{affinor} (\cite{9}%
) or \textit{vector $1$-form}; the remarkable affinor of every manifold is
the Kronecker tensor field $I=(\delta _{j}^{i})$.

\medskip

Fix $\mu \in \{-1, +1\}$. Let now $g$ be a Riemannian metric on $M$.

\medskip

\textbf{Definition 1.1} $M$ is a $(g, F, \mu )$-\textit{manifold} if : 
\begin{equation*}
g(FX, Y)+\mu g(X, FY)=0. \eqno(1.1)
\end{equation*}
The geometry of the data $(M, g, F, \mu )$ is called \textit{affinor-metric
geometry}. If in particular, $F_{x}$ is nondegenerate at any point $x\in M$
then we say that $M$ is a \textit{nondegenerate} $(g, F, \mu )$-\textit{%
manifold}; otherwise, $M$ is called \textit{degenerate} $(g, F, \mu )$-%
\textit{manifold}.

\medskip

The relation $(1.1)$ says that the $g$-adjoint of $F$ is $F^*=-\mu F$. In
literature there is an abundance of examples of $(g, F, \mu )$-manifolds.
Some of the main examples are presented here:

\medskip

\textbf{Examples 1.2} \newline
1. An \textit{almost Hermitian manifold} (\cite[p. 11]{4}) $(M,g,J)$ is a
nondegenerate $(g,F,\mu =+1)$-manifold; the nondegeneration is a consequence
of $J^{2}=-I$. \newline
2. An \textit{almost parahermitian manifold} (\cite{1}) $(M,g,P)$ is a
nondegenerate $(g,F,\mu =+1)$-manifold while an \textit{almost Riemannian
product manifold} is a nondegenerate $(g,F,\mu =-1)$-manifold; the
nondegeneration is a consequence of $P^{2}=I$. \newline
3. An \textit{almost contact metric manifold} (\cite[p. 15]{4}) $%
(M,g,\varphi ,\xi ,\eta )$ is a $(g,F,\mu =+1)$-manifold; as $\varphi (\xi
)=0$, $M$ is degenerate. \newline
4. An \textit{almost paracontact manifold} (\cite{12}) $(M,g,\varphi ,\xi
,\eta )$ is a $(g,F,\mu =+1)$-manifold. As in the previous example we have $%
\varphi (\xi )=0$ and therefore $M$ is degenerate.\newline
5. The general case of a nondegenerate $(g,F,\mu =+1)$-manifold is called 
\textit{structured manifold} in \cite{11}.

\medskip

Recall that a real $2m$-dimensional manifold $M$ is called an \textit{almost
symplectic manifold} if it is endowed with a nondegenerate $2$-form $\Omega
\in \Lambda ^2(M)$. We derive the following characterization:

\medskip

\textbf{Proposition 1.3} \textit{Let $M$ be a $(g, F, \mu =+1)$-manifold.
Then $M$ is nondegenerate if and only if $\Omega $ defined by: \newline
\begin{equation*}
\Omega (X, Y)=g(FX, Y) \eqno(1.2)
\end{equation*}
is an almost symplectic structure. In this case $m$ is even}.

\medskip

\textbf{Proof} $\Omega $ is skew-symmetric from $\mu =+1$. A straightforward
computation yields that $\Omega $ is nondegenerate if and only if $M$ is a
nondegenerate $(g, F, \mu =+1)$-manifold. \quad $\Box $

\medskip

\textbf{Example 1.4} For Example 1.2.1 $\Omega $ is exactly the \textit{%
fundamental} or \textit{K\"ahler} $2$-form and then inspired by this fact we
introduce:

\medskip

\textbf{Definition 1.5} For a nondegenerate $(g, F, \mu =+1)$-manifold $%
\Omega $ is call \textit{the fundamental} $2$-\textit{form}.

\medskip

In the last part of this section let us recall briefly the geometry of
Riemannian submanifolds. Consider an $n$-dimensional submanifold $N$ of $M$.
Then the main objects induced by the Levi-Civita connection $\widetilde{%
\nabla }$ of $(M,g)$ on $N$ are involved in the well known Gauss-Weingarten
equations: 
\begin{equation*}
\widetilde{\nabla }_{X}Y=\nabla _{X}Y+h(X,Y),\quad \widetilde{\nabla }%
_{X}V=-A_{V}X+\nabla _{X}^{\perp }V,\eqno(1.3)
\end{equation*}%
for any $X,Y\in \Gamma (TN)$ and $V\in \Gamma (T^{\bot }N)$. Here $\nabla $
is the Levi-Civita connection on $N$, $h$ is the second fundamental form of $%
N$, $A_{V}$ is the Weingarten operator with respect to the normal section $V$
and $\nabla ^{\bot }$ is the normal connection in the normal bundle $T^{\bot
}N$ of $N$. Let us point out that $h$ and $A_{V}$ are related by: 
\begin{equation*}
g(h(X,Y),V)=g(A_{V}X,Y).\eqno(1.4)
\end{equation*}%
If $h$ vanishes identically on $N$ then $N$ is called \textit{totally
geodesic}.

\section{Submanifolds in affinor-metric geometry}

Next, we consider a submanifold $N$ of a $(g,F,\mu )$-manifold $M$. Then $g$
induces a Riemannian metric on $N$ which we denote by the same symbol $g$.
Then, following the definition given by Bejancu \cite{2} for CR-submanifolds
we introduce a special class of submanifolds of $M$ as follows:

\medskip

\textbf{Definition 2.1} $N$ is a \textit{semi-invariant submanifold} of $M$
if there exists a distribution $D$ on $N$ satisfying the conditions: \newline
(i) $D$ is $F$-invariant: 
\begin{equation*}
F(D_{x})\subset D_{x}\ ,\ \ \ \forall \ x\in N. \eqno(2.1)
\end{equation*}
(ii) The complementary orthogonal distribution $D^{\perp }$ to $D$ in $TN$
is $F$-anti-invariant, that is: 
\begin{equation*}
F(D_{x}^{\perp })\subset T_{x}^{\bot}N,\ \ \ \forall \ x\in N. \eqno(2.2)
\end{equation*}
(iii) $F^{2}(D^{\perp })$ is a distribution on $N$. \newline
Some particular classes of semi-invariant submanifolds are defined as
follows. Let $p$ and $q$ be the ranks of the distributions $D$ and $D^{\perp
}$ respectively. If $q=0$, that is $D^{\perp }=\left\{ 0\right\} $, we say
that $N$ is an $F$-\textit{invariant submanifold} of $M$. If $p=0$, that is $%
D=\left\{ 0\right\} $, we call $N$ an $F$-\textit{anti-invariant submanifold}
of $M$.

If $pq\neq 0$ then $N$ is called a \textit{proper} semi-invariant
submanifold. Now, we denote by $\widetilde{D}$ the complementary orthogonal
vector bundle to $F(D^{\perp })$ in $T^{\perp }N$. If $\widetilde{D}%
=\left\{0\right\} $ then we say that $N$ is a \textit{normal} $F$%
-semi-invariant submanifold.

\smallskip

Thus, $N$ is an $F$-invariant, respectively $F$-anti-invariant, if and only
if: 
\begin{equation*}
F(TN)\subset TN\ \ \ (\text{resp. }F(TN)\subset T^{\perp }N). \eqno(2.3)
\end{equation*}
$N$ is normal $F$-semi-invariant if and only if: 
\begin{equation*}
F(D^{\perp })=T^{\perp }N. \eqno(2.4)
\end{equation*}

\smallskip

\textbf{Examples 2.2} \newline
1) For Example 1.2.1 we obtain the classical concept of CR-submanifold of
Bejancu \cite[p. 20]{4}; the condition iii) is satisfied from $J^{2}=-I$.%
\newline
2) For Example 1.2.2 we obtain the notion of semi-invariant submanifold; for
the almost parahermitian case see \cite{1} while for the second case see 
\cite{3}. The condition iii) is satisfied again from $P^{2}=-I$.\newline
3) For Example 1.2.3 we obtain the notion of semi-submanifold \cite[p. 100]%
{4} with $\xi \in T^{\bot }N$. This last condition implies $TN\subset ker\
\eta $ and since $\varphi |_{ker\ \eta }$ is an almost complex structure we
get iii).\newline
4) For Example 1.2.4 we obtain the concept of semi-submanifold from \cite{10}
with $\xi \in T^{\bot }N$. Again this condition means $TN\subset ker\ \eta $
and since $\varphi |_{ker\ \eta }$ is an almost product structure we have
iii). \newline
5) The condition (iii) does not appears in \cite{11}.

\medskip

Returning to the Definition 2.1 we deduce that the tangent bundle $TN$ and
the normal bundle $T^{\bot }N$ of a semi-invariant submanifold $N$ have the
orthogonal decompositions: 
\begin{equation*}
TN=D\oplus D^{\perp }, \quad T^{\perp }N=F(D^{\perp })\oplus \widetilde{D}. %
\eqno(2.5)
\end{equation*}
Then we denote by $P$ and $Q$ the projection morphisms of $TN$ on $D$ and $%
D^{\perp }$ respectively and obtain for $X=PX+QX\in \Gamma (TN)$: 
\begin{equation*}
FX=\varphi X+\omega X \eqno(2.6)
\end{equation*}
where we put: 
\begin{equation*}
\varphi =F\circ P, \quad \omega =F\circ Q. \eqno(2.7)
\end{equation*}
Thus $\varphi $ is a tensor field of $(1, 1)$-type on $N$ while $\omega $ is
a $F(D^{\perp })-$valued vector $1$-form on $N$. Thus we derive:

\medskip

\textbf{Proposition 2.5} \textit{Let} $N$ \textit{be a semi-invariant
submanifold of a} $(g, F, \mu )$-\textit{manifold} $M$. \textit{Then}: 
\newline
(iv) $N$ \textit{is a} $(g, \varphi , \mu )$-\textit{manifold}. \newline
(v) $F^{2}(D^{\perp })$ \textit{is a vector subbundle of} $D^{\perp }$. 
\newline
(vi) \textit{The vector bundle} $\widetilde{D}$ \textit{is} $F$-\textit{%
invariant i.e. for all} $x\in N$ \textit{we have}: $F(\widetilde{D}%
_{x})\subset \widetilde{D}_{x}$.

\medskip

\textbf{Proof} (iv) By definition, $g$ is a Riemannian metric on $N$ and $%
\varphi $ is a tensor field of $(1,1)$-type on $N$; we need only to show $%
(2.1)$. By using $(1.1)$ for $F$ we obtain for $X,Y\in \Gamma (TN)$: 
\begin{eqnarray*}
g(\varphi X,Y) &=&g(FPX,Y)=g(FPX,PY)=-\mu g(PX,FPY)= \\
&=&-\mu g(X,FPY)=-\mu g(X,\varphi Y).
\end{eqnarray*}%
(v) Take $X\in \Gamma (D)$ and $Y\in \Gamma (D^{\perp })$ in $(2.1)$: $%
g(X,F^{2}Y)=-\mu g(FX,FY)=0$ since $FX\in \Gamma (D)$ and $FY\in \Gamma
(T^{\perp }N)$. Hence $F^{2}(D^{\perp })$ is orthogonal to $D$ and by
condition (iii) we deduce that $F^{2}(D^{\perp })$ is a vector subbundle of $%
D^{\perp }$. \newline
(vi) Take $X\in \Gamma (TN)$, $Y\in \Gamma (D^{\perp })$ and $V\in \Gamma (%
\widetilde{D})$. Then we obtain: 
\begin{equation*}
g(FV,X)=-\mu g(V,FX)=-\mu g(V,\varphi X+\omega X)=0
\end{equation*}%
and: 
\begin{equation*}
g(FV,FY)=-\mu g(V,F^{2}Y)=0
\end{equation*}%
since $\varphi X\in \Gamma (D)$, $\omega X\in \Gamma (FD^{\perp })$ and $%
F^{2}Y\in \Gamma (D^{\perp })$. Thus $F\widetilde{D}$ is orthogonal to $%
TN\oplus FD^{\perp }$, that is $F\widetilde{D}$ is a vector subbundle of $%
\widetilde{D}$. This completes the proof of the proposition. \quad $\Box $

\medskip

In the non-degenerated case we have equalities for the above inclusions:

\medskip

\textbf{Corollary 2.6} \textit{Let} $N$ \textit{be a semi-invariant
submanifold of a nondegenerate} $(g, F, \mu )$-\textit{manifold} $M$. 
\textit{Then}: \newline
1) \textit{the above distributions satisfy}: 
\begin{equation*}
F(D)=D, \quad F^{2}(D^{\perp })=D^{\perp }, \quad F(\widetilde{D})=%
\widetilde{D}. \eqno(2.8)
\end{equation*}
2) \textit{if $\mu =+1$ then $D^{\bot }$ and $F(D^{\bot })$ are Lagrangian
distribution on $(TM, \Omega )$. In particular if $N$ is a normal
semi-invariant submanifold then $T^{\bot }N$ is a Lagrangian submanifold in $%
(TM, \Omega )$}.

\medskip

\textbf{Proof} We need to prove only 2). \newline
2.1) Let $X,Y\in \Gamma (D^{\bot })$; then $\Omega (X,Y)=g(FX,Y)=0$ since $%
FX\in \Gamma (T^{\bot }N)$ while $Y\in \Gamma (TN)$. \newline
2.2) Let $X,Y\in \Gamma (F(D^{\bot }))$; then $\Omega (X,Y)=g(FX,Y)=0$ since 
$FX\in \Gamma (TN)$ while $Y\in \Gamma (T^{\bot }N)$. \quad $\Box $

\medskip

The second part of the above Corollary is extremely important since it
relates the geometry of semi-invariant submanifolds with the almost
symplectic geometry, a topic very studied from the point of view of
applications in Analytical Mechanics.

\section{Integrability of distributions on a semi-invariant submanifold}

Let $N$ be a semi-invariant submanifold of a $(g,F,\mu )$-manifold $M$. Then
we recall that the Nijenhuis tensor field of $F$ is defined as follows (\cite%
[p.11]{4}): 
\begin{equation*}
N_{F}(X,Y)=[FX,FY]+F^{2}[X,Y]-F[FX,Y]-F[X,FY],\eqno(3.1)
\end{equation*}%
for any $X,Y\in \Gamma (TM)$. In a similar way, the Nijenhuis tensor field
of $\varphi $ on $N$ is given by: 
\begin{equation*}
N_{\varphi }(X,Y)=[\varphi X,\varphi Y]+\varphi ^{2}[X,Y]-\varphi \lbrack
\varphi X,Y]-\varphi \lbrack X,\varphi Y],\eqno(3.2)
\end{equation*}%
for any $X,Y\in \Gamma (TN)$. We recall that a tensor field of $(1,1)$-type
defines an \textit{integrable structure} on a manifold if and only if its
Nijenhuis tensor field vanishes identically on the manifold. Now we obtain
necessary and sufficient conditions for the integrability of $D$ and $%
D^{\perp }$ in terms of Nijenhuis tensor fields of $F$ and $\varphi $.

\medskip

\textbf{Theorem 3.1} \textit{Let} $N$ \textit{be a semi-invariant
submanifold of a} $(g, F, \mu )$-\textit{manifold} $M$. \textit{Then the
following assertions are equivalent}: \newline
1) $D$ \textit{is an integrable distribution}. \newline
2) \textit{The Nijenhuis tensor field of} $\varphi $ \textit{satifies}: 
\begin{equation*}
Q\circ N_{\varphi }=0, \quad \forall X, Y\in \Gamma(D). \eqno(3.3)
\end{equation*}
3) \textit{The Nijenhuis tensor fields of} $F$ \textit{and} $\varphi $ 
\textit{satisfy the equality}: $N_F=N_{\varphi }$ \textit{on} $D$.

\medskip

\textbf{Proof} Firstly, we note that $D$ is integrable if and only if: 
\begin{equation*}
Q([X,Y])=0,\quad \forall X,Y\in \Gamma (D).\eqno(3.4)
\end{equation*}%
Since the last three terms in the right side of $(3.2)$ lie in $\Gamma (D)$
we deduce that: 
\begin{equation*}
Q\circ N_{\varphi }(X,Y)=Q([FX,FY]),\quad \forall X,Y\in \Gamma (D).\eqno%
(3.5)
\end{equation*}%
As $M$ is nondegenerate we deduce that $\varphi $ is an automorphism on $%
\Gamma (D)$. Thus the equivalence of 1) and 2) follows directly. Next, we
obtain for any $X,Y\in \Gamma (D)$: 
\begin{equation*}
N_{F}(X,Y)=N_{\varphi }(X,Y)+F\omega ([X,Y])-\omega ([\varphi X,Y])-\omega
([X,\varphi Y]).\eqno(3.6)
\end{equation*}%
If $D$ is integrable then the last three terms of $(3.6)$ vanishes and this
yields 3). Conversely, suppose that $N_{F}=N_{\varphi }$ on $D$; then: 
\begin{equation*}
F\omega ([X,Y])=\omega ([\varphi X,Y]+[X,\varphi Y]).\eqno(3.7)
\end{equation*}%
Obviously the right-hand-side of the previous equation is in $\Gamma
(F(D^{\bot }))\subset \Gamma (T^{bot}N)$. On the other hand, the
left-hand-side is in $\Gamma (F^{2}D^{\bot })\subset \Gamma (TN)$; we
conclude that both sides in $(3.7)$ must vanish. 

Finally, from:\noindent {} $F^{2}Q([X,Y])=0$ and $F^{2}$ automorphism of $%
\Gamma (TM)$ we deduce 1). \quad $\ \ \Box $

\medskip

\textbf{Remark 3.2} For Example 1.2.1 the equivalence of 1) and 2) is
exactly the Theorem 2.2. of \cite[p. 25]{4} while the equivalence of 1) and
3) is the Theorem 2.1. of \cite[p. 25]{4}.

\medskip

Now, we consider $X, Y\in \Gamma (D^{\bot })$. Then taking into account that 
$\varphi X=\varphi Y=0$ we get: 
\begin{equation*}
N_{\varphi }(X, Y)=F^{2}P[X, Y] \eqno(3.8)
\end{equation*}
and this enables us to state the following:

\medskip

\textbf{Theorem 3.3} \textit{Let} $N$ \textit{be a semi-invariant
submanifold of a nondegenerate} $(g, F, \mu )-$\textit{manifold. Then} $%
D^{\bot }$ \textit{is integrable if and only if the Nijenhuis tensor field of%
} $\varphi $ \textit{vanishes identically on} $D^{\bot }$.

\medskip

\textbf{Remark 3.4} For Example 1.2.1 the above result is the Theorem 2.3.
of \cite[p. 26]{4}.

\section{A natural foliation on a semi-invariant submanifold}

Let $\widetilde{\nabla }$ be the Levi-Civita connection on $M$ with respect
to the Riemannian metric $g$. Then $F$ is a \textit{parallel tensor field}
on $M$ if: 
\begin{equation*}
\widetilde{\nabla }F=0. \eqno(4.1)
\end{equation*}

\smallskip

\textbf{Examples 4.1}\newline
1) For Example 1.2.1 we have the notion of \textit{K\"{a}hler manifold}. 
\newline
2) For Example 1.2.2, in the first part we have the concept of \textit{para-K%
\"{a}hler manifold} while for the second part the notion of \textit{locally
Riemannian product manifold}. \newline
3) For Example 1.2.3 we get the notion of \textit{cosymplectic manifold}.

\medskip

In the present section we study the geometry of semi-invariant submanifolds
of $(g, F, \mu )$-manifolds with parallel tensor field $F$. First, we prove
the following:

\medskip

\textbf{Proposition 4.2} \textit{Let} $N$ \textit{be a semi-invariant
submanifold of a nondegenerate} $(g, F, \mu )$-\textit{manifold with
parallel tensor field} $F$. \textit{Then for all} $X, Y\in \Gamma (D^{\bot
}) $: 
\begin{equation*}
A_{FX}Y-A_{FY}X=\varphi ([X, Y]). \eqno(4.2)
\end{equation*}

\smallskip

\textbf{Proof} By using the Weingarten equation and the parallelism
condition we get: 
\begin{equation*}
A_{FX}Y=\nabla ^{\bot}_{Y}FX-\nabla _{Y}FX=\nabla ^{\bot}_{Y}FX-F(\widetilde{%
\nabla }_XY). \eqno(4.2)
\end{equation*}
Writing a similar equation by interchanging $X$ and $Y$ and then subtracting
we obtain: 
\begin{equation*}
A_{FX}Y-A_{FY}X=\nabla _{Y}^{\bot }FX-\nabla _{X}^{\bot }FY+F([X, Y]), \eqno%
(4.3)
\end{equation*}
since $\nabla $ is a torsion-free linear connection. Thus $(4.2)$ is
obtained by equalizing the tangent parts to $N$ in the above equation. \quad 
$\Box $

\medskip

\textbf{Example 4.3} The relation $(4.2)$ becomes for Example 1.2.1 the
equation $(2.2)$ of \cite[p. 43]{4}.

\medskip

Now, we can state the following main result:

\medskip

\textbf{Theorem 4.4} \textit{Let} $N$ \textit{be a semi-invariant
submanifold of a nondegenerate} $(g, F, \mu =+1)$-\textit{manifold with
parallel tensor field} $F$. \textit{Then the} $F$-\textit{anti-invariant
distribution} $D^{\bot }$ \textit{is integrable}.

\medskip

\textbf{Proof} For any $X, Y\in \Gamma (D^{\bot })$ and $Z\in \Gamma (D)$ we
have: 
\begin{equation*}
g(A_{FX}Y, Z)=-g(F\widetilde{\nabla }_{Y}X, Z)=+\mu g(\widetilde{\nabla }%
_{Y}X, FZ)=-\mu g(X, \widetilde{\nabla }_{Y}FZ)=
\end{equation*}
\begin{equation*}
=\mu ^2g(FX, \widetilde{\nabla }_{Y}Z)=\mu ^2g(FX, [Y,Z]+\widetilde{\nabla }%
_{Z}Y)=\mu ^2g(FX, \widetilde{\nabla }_{Z}Y). \eqno(4.4)
\end{equation*}
Also, we have: 
\begin{equation*}
g(A_{FY}X, Z)=\mu ^2g(FY, \widetilde{\nabla }_{Z}X)=-\mu ^2g(F\widetilde{%
\nabla }_{Z}Y, X)=\mu ^3g(\widetilde{\nabla }_{Z}Y, FX). \eqno(4.5)
\end{equation*}
Comparing $(4.4)$ and $(4.5)$ we deduce that for $\mu =+1$: 
\begin{equation*}
g(A_{FX}Y-A_{FY}X, Z)=0
\end{equation*}
which means that $A_{FX}Y-A_{FY}X\in \Gamma (D^{\bot })$. On the other hand,
from $(4.2)$ we conclude that: 
\begin{equation*}
A_{FX}Y-A_{FY}X\in \Gamma (D).
\end{equation*}
and thus we have that: 
\begin{equation*}
A_{FX}Y-A_{FY}X=0. \eqno(4.6)
\end{equation*}
Finally, returning to $(4.2)$ and taking into account that $F$ is
nondegenerate we deduce that: 
\begin{equation*}
P[X, Y]=0,
\end{equation*}
that is, $D^{\bot }$ is integrable. \quad $\Box $

\medskip

\textbf{Remark 4.5} For Example 1.2.1 the above result is part (i) of
Theorem 1.1. of \cite[p. 39]{4}.

\medskip

Regarding the integrability of $D$ we prove the following:

\medskip

\textbf{Theorem 4.6} \textit{Let} $N$ \textit{be a semi-invariant
submanifold of a nondegenerate} $(g, F, \mu )$-\textit{manifold} $M$\textit{%
with parallel tensor field} $F$. \textit{Then the} $F$-\textit{invariant
distribution} $D$ \textit{is integrable if and only if the second
fundamental form} $h$ \textit{of} $N$ \textit{satisfies for any} $X, Y\in
\Gamma (D)$ \textit{and} $Z\in \Gamma (D^{\bot })$: 
\begin{equation*}
g(h(X, \varphi Y)-h(Y, \varphi X), FZ)=0. \eqno(4.7)
\end{equation*}

\smallskip

\textbf{Proof} By using the Gauss equation we deduce that: 
\begin{equation*}
\nabla _{X}\varphi Y+h(X, \varphi Y)=\varphi (\nabla _{X}Y)+\omega (\nabla
_{X}Y)+Fh(X,Y). \eqno(4.8)
\end{equation*}
Write a similar equation by interchanging $X$ and $Y$, and then subtracting
we obtain: 
\begin{equation*}
\nabla _{X}\varphi Y-\nabla _{Y}\varphi X+h(X, \varphi Y)-h(Y, \varphi
X)=\varphi ([X, Y])+\omega ([X,Y]) \eqno(4.9)
\end{equation*}
since $h$ is symmetric and $\nabla $ is a torsion-free linear connection.
Equalize the normal parts in the above equation and obtain: 
\begin{equation*}
h(X, \varphi Y)-h(Y, \varphi X)=\omega ([X, Y]). \eqno(4.10)
\end{equation*}
Now, suppose that $D$ is integrable; then $(4.7)$ is immediately.
Conversely, if $(4.6)$ is satisfied, then with $(4.10)$ we deduce that: 
\begin{equation*}
0=-\mu g(Q[X, Y], F^{2}Z).
\end{equation*}
Since $F$ is nondegenerate we infer that $F^{2}$ is an automorphism of $%
\Gamma (D^{\bot })$ and hence $D$ is integrable. \quad $\Box $

\medskip

\textbf{Remark 4.7} In particular, if $F$ is an almost complex structure on $%
M$ then we obtain the results of Bejancu \cite{4} and Blair-Chen \cite{6}
respectively, for CR-submanifolds.

\medskip

Now, for $\mu =+1$ we denote by $\mathcal{F}^{\bot }$ the natural foliation
defined by the $F$-anti-invariant distribution $D^{\bot }$ and call it the $%
F $-\textit{anti-invariant foliation} on $N$. We recall that $\mathcal{F}%
^{\bot }$ is called a \textit{totally geodesic foliation} if each leaf of $%
\mathcal{F}^{\bot }$ is totally geodesic immersed in $N$. Thus $\mathcal{F}%
^{\bot }$ is totally geodesic if and only if the Levi-Civita connection $%
\nabla $ of $N$ satisfies for all $Y, Z\in \Gamma (D^{\bot })$: 
\begin{equation*}
\nabla _{Y}Z\in \Gamma (D^{\bot }). \eqno(4.11)
\end{equation*}

\smallskip

\textbf{Theorem 4.8} \textit{Let} $N$ \textit{be a semi-invariant
submanifold of a nondegenerate} $(g, F, \mu )$-\textit{manifold} $M$ \textit{%
with parallel tensor field} $F$. \textit{Then the following assertions are
equivalent}: \newline
(i) \textit{The} $F$-\textit{anti-invariant foliation is totally geodesic}. 
\newline
(ii) \textit{The second fundamental form} $h$ \textit{of} $N$ \textit{%
satisfies for all} $X\in \Gamma (D)$ \textit{and} $Y\in \Gamma (D^{\perp })$%
: 
\begin{equation*}
h(X, Y)\in \Gamma (\widetilde{D}). \eqno(4.12)
\end{equation*}
(iii) $D^{\bot }$ \textit{is} $A_{V}$-\textit{invariant for any} $V\in
\Gamma (FD^{\bot })$ \textit{that is we have for all} $Y\in \Gamma (D^{\bot
})$: 
\begin{equation*}
A_{V}Y\in \Gamma (D^{\bot }). \eqno(4.13)
\end{equation*}

\smallskip

\textbf{Proof} We have for any $X\in \Gamma (D)$ and $Y,Z\in \Gamma (D^{\bot
})$: 
\begin{equation*}
g(\nabla _{Y}Z,FX)=g(\widetilde{\nabla }_{Y}Z,FX)=-\mu g(\widetilde{\nabla }%
_{Y}FZ,X)=
\end{equation*}

\begin{equation*}
=\mu g(A_{FZ}Y,X)=\mu g(h(X,Y),FZ).\eqno(4.14)
\end{equation*}%
Now, suppose that $\mathcal{F}^{\bot }$ is totally geodesic; then the first
term of $(4.14)$ vanishes. Hence the last term in $(4.14)$ vanishes which
implies ii). Conversely, suppose $(4.12)$ is satisfied. Then from $(4.14)$
we deduce $(4.11)$ since $F$ is an automorphism of $\Gamma (D)$. This proves
the equivalence of (i) and (ii). The equivalence of (ii) and (iii) is
straightforward. \quad $\Box $

\medskip

\textbf{Remark 4.9} For Example 1.2.1 the equivalence of (i) and (ii) is the
Theorem 1.3. of \cite[p. 41]{4}.

\medskip

Finally, we can prove the following:

\medskip

\textbf{Theorem 4.10} \textit{Let} $N$ \textit{be a semi-invariant
submanifold of a nondegenerate} $(g, F, \mu )$-\textit{manifold with
parallel tensor field} $F$. \textit{Then the} $F$-\textit{invariant
distribution} $D$ \textit{is integrable and the foliation} $\mathcal{F}$ 
\textit{defined by} $D$ \textit{is totally geodesic if and only if the
second fundamental form} $h$ \textit{of} $N$ \textit{satisfies for all} $X,
Y\in \Gamma (D)$: 
\begin{equation*}
h(X, Y)\in \Gamma (\widetilde{D}). \eqno(4.15)
\end{equation*}

\smallskip

\textbf{Proof} $D$ is integrable and $\mathcal{F}$ is totally geodesic if
and only if for all $X,U\in \Gamma (D)$: 
\begin{equation*}
\nabla _{X}U\in \Gamma (D).\eqno(4.16)
\end{equation*}%
This is equivalent to: 
\begin{equation*}
g(\widetilde{\nabla }_{X}U,Z)=0,\eqno(4.17)
\end{equation*}%
for all $Z\in \Gamma (D^{\bot })$ As $F$ is an automorphism of $\Gamma (D)$
we can write the above equality as follows: 
\begin{equation*}
g(\widetilde{\nabla }_{X}FY,Z)=0,\eqno(4.18)
\end{equation*}%
for all $X,Y\in \Gamma (D)$ and $Z\in \Gamma (D^{\bot })$, which is
equivalent to: 
\begin{equation*}
g(\widetilde{\nabla }_{X}Y,FZ)=0.\eqno(4.19)
\end{equation*}%
By using the Gauss equation, the last relation is equivalent to: 
\begin{equation*}
g(h(X,Y),FZ)=0,\eqno(4.20)
\end{equation*}%
which completes the proof of the theorem. \quad $\Box $

\medskip

\textbf{Remark 4.11} For Example 1.2.1 the above result is Theorem 1.2. of 
\cite[p. 40]{4}.

\end{document}